\documentclass[11pt]{article}

\usepackage{amsmath,amsthm,amsfonts,amssymb,latexsym}
\usepackage{amsmath}

\usepackage{a4wide}

\newtheorem{thm}{Theorem}[section]
\newtheorem{cor}[thm]{Corollary}
\newtheorem{lem}[thm]{Lemma}
\newtheorem{prop}[thm]{Proposition}

\newtheorem{theorem}{Theorem}[section]

\newtheorem{proposition}[theorem]{Proposition}
\newtheorem{corollary}[theorem]{Corollary}

\newcommand{\bbR}{\mathbb{R}}

\newcommand{\fF}{{\mathbb{F}}}

\newcommand{\fL}{{\mathbb{L}}}

\newcommand{\cK}{{\mathcal{K}}}

\newcommand{\cF}{{\mathcal{F}}}

\newcommand{\Hom}{\mathrm{Hom}}

\newcommand{\HH}{\mathrm{H}}

\newcommand{\half}{\frac{1}{2}}

\def\a{\alpha}
\def\b{\beta}

\def\l{\lambda}

\begin{document}


\title{The Binary Invariant Differential Operators on Weighted Densities on
the superspace $\mathbb{R}^{1|n}$ and Cohomology }


\author{Mabrouk Ben Ammar\thanks{D\'epartement de Math\'ematiques,
Facult\'e des Sciences de Sfax, BP 802, 3038 Sfax, Tunisie.
E. mail: mabrouk.benammar@fss.rnu.tn,}\and    Nizar Ben Fraj\thanks{Institut
Sup\'{e}rieur de Sciences Appliqu\'{e}es et Technologie, Sousse, Tunisie. E. mail: benfraj\_nizar@yahoo.fr}\and Salem Omri\thanks{D\'epartement de Math\'ematiques, Facult\'e des Sciences de
Gafsa, Zarroug 2112 Gafsa, Tunisie. E. mail: omri\_salem@yahoo.fr,}}

\maketitle

\begin{abstract}
Over the $(1,n)$-dimensional real superspace, $n>1$, we classify
 $\mathcal{K}(n)$-invariant binary differential operators acting on
the superspaces of weighted densities, where $\mathcal{K}(n)$ is
the Lie superalgebra of contact vector fields. This result allows
us to compute the first differential cohomology of
$\mathcal{K}(n)$ with coefficients in the superspace
of linear differential operators acting on the superspaces of
weighted densities--a superisation of a result by Feigin and
Fuchs. We explicitly give 1-cocycles spanning these cohomology
spaces.
\end{abstract}

\maketitle {\bf Mathematics Subject Classification} (2000). 53D55




\section{Introduction}
This work is a direct continuation of \cite{DYA, GP} and \cite{N,
c} listed among other things, binary differential operators
invariant with respect to a supergroup of diffeomorphisms and
computed cohomology of polynomial versions of various infinite
dimensional Lie superalgebras.

Let $\mathfrak{vect}(1)$ be the Lie algebra of polynomial vector
fields on $\mathbb{R}$. Consider the 1-parameter deformation of
the $\mathfrak{vect}(1)$-action  on $\mathbb{R}[x]$:
\begin{equation*}
L_{X\frac{d}{dx}}^\lambda(f)= Xf'+\lambda X'f,
\end{equation*}
where $X, f\in\mathbb{R}[x]$ and $X':=\frac{dX}{dx}$. Denote by
$\cF_\l$ the $\mathfrak{vect}(1)$-module structure on $\bbR[x]$
defined by $L^\lambda$ for a fixed $\lambda$. Geometrically,
${\cal F}_\lambda=\left\{fdx^{\lambda}\mid f\in
\mathbb{R}[x]\right\}$ is the space of polynomial weighted
densities of weight $\l\in\mathbb{R}$. The space ${\cal
F}_\lambda$ coincides with the space of vector fields, functions
and differential 1-forms for $\lambda = -1,\, 0$ and $1$,
respectively.

Denote by $\mathrm{D}_{\lambda,\mu}:=\Hom_{\rm{diff}}({\cal F}_\lambda,
{\cal F}_\mu)$ the $\mathfrak{vect}(1)$-module of linear
differential operators with the natural
$\mathfrak{vect}(1)$-action. Feigin and Fuchs \cite{ff} computed
$\HH^1_{\rm diff}\left(\mathfrak{vect}(1);
\mathrm{D}_{\lambda,\mu}\right)$, where $\mathrm{H}^*_\mathrm{diff}$
denotes the differential cohomology; that is, only cochains given
by differential operators are considered. They showed that
non-zero cohomology $\mathrm{H}^1_{\rm
diff}\left(\mathfrak{vect}(1);\mathrm{D}_{\lambda,\mu}\right)$ only
appear for particular values of weights that we call {\it
resonant} which satisfy $\mu-\lambda\in\mathbb{N}$. These spaces arise
in the classification of infinitesimal deformations of the
$\mathfrak{vect}(1)$-module ${\cal
S}_{\mu-\lambda}=\bigoplus_{k=0}^\infty{\cal F}_{\mu-\lambda-k}$, the
space of symbols of $\mathrm{D}_{\lambda,\mu}$.

On the other hand, Grozman \cite{GP} classified all $\mathfrak{vect}(1)$-invariant
binary differential operators on $\mathbb{R}$ acting in the spaces
$\mathcal{F}_\lambda$. He showed that all
invariant operators are of order $\leq3$ and can be expressed as a
composition of the Rham differential and the Poisson bracket,
except for one called Grozman operator.

It is natural to study the simplest super analog of the problems
solved respectively in [5] and [10], namely, we consider the
superspace $\mathbb{R}^{1|n}$ endowed with its standard contact
structure defined by the 1-form $\a_n$, and the Lie superalgebra
$\cK(n)$ of contact polynomial vector fields on
$\mathbb{R}^{1|n}$. We introduce the $\cK(n)$-module
$\mathbb{F}_\lambda^n$ of $\lambda$-densities on
$\mathbb{R}^{1|n}$ and the $\cK(n)$-module of linear differential
operators, $\mathbb{D}^n_{\lambda,\mu}
:=\Hom_{\rm{diff}}(\fF_{\lambda}^n,\fF_{\mu}^n)$, which are super
analogues of the spaces $\mathcal{F}_\lambda$ and
$\mathrm{D}_{\lambda,\mu}$, respectively. The classification of
the $\cK(1)$-invariant binary differential operators on
$\mathbb{R}^{1|1}$ acting in the spaces $\mathbb{F}_\lambda^1$ is
due to Leites et al. \cite{DYA}, while the space $\HH^1_{\rm
diff}\left(\cK(1);\mathbb{D}^1_{\lambda,\mu}\right)$ has been
computed by Basdouri et al. \cite{bbbbk} (see also \cite{c}) and
the space $\HH^1_{\rm
diff}\left(\cK(2);\mathbb{D}^2_{\lambda,\mu}\right)$ has been
computed by the second author \cite{N}. We also mention that Duval
and Michel studied a similar problem for the group of
contactomorphisms of the supercircle $S^{1\mid n}$ instead of
$\cK(n)$ related to the link between discrete projective
invariants of the supercircle, and the cohomology of the group of
its contactomorphisms \cite{D}.

In this paper we classify all
$\cK(n)$-invariant binary differential operators on
$\mathbb{R}^{1|n}$ acting in the spaces
$\mathbb{F}_\lambda^n$ for $n>1$. We use
the result to compute $\HH^1_{\rm
diff}\left(\cK(n);\mathbb{D}^n_{\lambda,\mu}\right)$ for $n>2$. We
show that, as in the classical setting, non-zero cohomology
$\HH^1_{\rm diff}\left(\cK(n);\mathbb{D}^n_{\lambda,\mu}\right)$ only
appear for resonant values of weights which satisfy
$\mu-\lambda\in\frac{1}{2}\mathbb{N}$. Moreover, we give explicit basis
of these cohomology spaces. These spaces arise in the
classification of infinitesimal deformations of the
$\cK(n)$-module $ {\frak S}^n_{\mu-\lambda}=\bigoplus_{k\geq0}{\frak
F}^n_{\mu-\lambda-\frac{k}{2}}$, a super analogue of ${\cal
S}_{\mu-\lambda}$, see \cite{bbbbk}.

\section{Definitions and Notation}

\subsection{The Lie superalgebra of contact vector fields on
$\mathbb{R}^{1|n}$}

Let $\mathbb{R}^{1\mid n}$ be the superspace with coordinates
$(x,~\theta_1,\ldots,\theta_n),$ where
$\theta_1,\,\dots,\,\theta_n$ are the odd variables, equipped with
the standard contact structure given by the following $1$-form:
\begin{equation}
\label{a} \a_n=dx+\sum_{i=1}^n\theta_id\theta_i.
\end{equation}
On the space
$\mathbb{R}[x,\theta]:=\mathbb{R}[x,\theta_1,\dots,\theta_n]$, we
consider the contact bracket
\begin{equation}
\{F,G\}=FG'-F'G-\frac{1}{2}(-1)^{|F|}\sum_{i=1}^n{\eta}_i(F)\cdot
{\eta}_i(G),
\end{equation}where
${\eta}_i=\frac{\partial}{\partial
{\theta_i}}-\theta_i\frac{\partial}{\partial x}$ and $|F|$ is the
parity of $F$. Note that the derivations $\eta_i$ are the
generators of n-extended supersymmetry and generate the kernel of
the form (\ref{a}) as a module over the ring of polynomial
functions. Let $\mathrm{Vect_{Pol}}(\mathbb{R}^{1|n})$ be the
superspace of polynomial vector fields on ${\mathbb{R}}^{1|n}$:
\begin{equation*}\mathrm{Vect_{Pol}}(\mathbb{R}^{1|n})=\left\{F_0\partial_x
+ \sum_{i=1}^n F_i\partial_i \mid ~F_i\in\mathbb{R}[x,\theta]~
\text{ for all } i  \right\},\end{equation*} where
$\partial_i=\frac{\partial}{\partial\theta_i}$ and
$\partial_x=\frac{\partial}{\partial x} $, and consider the
superspace $\mathcal{K}(n)$ of contact polynomial vector fields on
${\mathbb{R}}^{1|n}$. That is, $\mathcal{K}(n)$ is the superspace
of vector fields on $\mathbb{R}^{1|n}$ preserving the distribution
singled out by the $1$-form $\alpha_n$: $$
\mathcal{K}(n)=\big\{X\in\mathrm{Vect_{Pol}}(\mathbb{R}^{1|n})~|~\hbox{there
exists}~F\in {\mathbb{R}}[x,\,\theta]~ \hbox{such
that}~{L}_X(\alpha_n)=F\alpha_n\big\}. $$ The Lie superalgebra
$\mathcal{K}(n)$ is spanned by the fields of the form:
\begin{equation*}
X_F=F\partial_x
-\frac{1}{2}\sum_{i=1}^n(-1)^{|F|}{\eta}_i(F){\eta}_i,\;\text{where
$F\in \mathbb{R}[x,\theta]$.}
\end{equation*}
In particular, we have $\cK(0)=\mathfrak{vect}(1)$. Observe that
${L}_{X_F}(\alpha_n)=X_1(F)\alpha_n$. The bracket in
$\mathcal{K}(n)$ can be written as:
\begin{equation*}
[X_F,\,X_G]=X_{\{F,\,G\}}.
\end{equation*}

\subsection{Modules of weighted densities}
We introduce a one-parameter family of modules over the Lie
superalgebra $\mathcal{K}(n)$. As vector spaces all these modules
are isomorphic to ${\mathbb{R}}[x,\,\theta]$, but not as
$\mathcal{K}(n)$-modules.

For every contact polynomial vector field $X_F$, define a
one-parameter family of first-order differential operators on
$\mathbb{R}[x,\theta]$:
\begin{equation}
\label{superaction} \fL^{\lambda}_{X_F}=X_F + \lambda F',\quad
\l\in\mathbb{R}.
\end{equation}
We easily check that
\begin{equation}
\label{crochet}
[\fL^{\lambda}_{X_F},\fL^{\lambda}_{X_G}]=\fL^\l_{X_{\{F,G\}}}.
\end{equation}
We thus obtain a one-parameter family of $\cK(n)$-modules on
$\mathbb{R}[x,\theta]$ that we denote $\mathbb{F}^n_\lambda$, the
space of all polynomial weighted densities on $\mathbb{R}^{1|n}$
of weight $\l$ with respect to $\a_n$:
\begin{equation}
\label{densities} \fF^n_\l=\left\{F\a_n^{\l} \mid F
\in\mathbb{R}[x,\theta]\right\}.
\end{equation}
In particular, we have $\mathbb
{F}_{\lambda}^0=\mathcal{F}_\lambda$. Obviously the adjoint
$\cK(n)$-module is isomorphic to the space of weighted densities
on $\mathbb{R}^{1|n}$ of weight $-1.$

\subsection{Differential operators on weighted densities}

A differential operator on $\mathbb{R}^{1|n}$ is an operator on
$\mathbb{R}[x,\theta]$ of the form:
\begin{equation}\label{diff}
A=\sum_{k=0}^M
\sum_{\varepsilon=(\varepsilon_1,\cdots,\varepsilon_n)}a_{k,\epsilon}(x,\theta)\partial_x^k
\partial_{1}^{\varepsilon_1}\cdots\partial_{n}^{\varepsilon_n};\,\,
\varepsilon_i=0,1;\,\,M\in\mathbb{N}.
\end{equation}
Of course any differential operator defines a linear mapping
$F\alpha_n^\lambda\mapsto(AF)\alpha_n^\mu$ from
$\mathbb{F}^n_{\lambda}$ to $\mathbb{F}^n_{\mu}$ for any
$\lambda,\,\mu\in\mathbb{R}$, thus the space of differential
operators becomes a family of ${\mathcal K}(n)$-modules $
\mathbb{D}^n_{\lambda,\mu}$ for the natural action:
\begin{equation}\label{d-action}
{X_F}\cdot A=\mathbb{L}^{\mu}_{X_F}\circ A-(-1)^{|A||F|} A\circ
\mathbb{L}^{\lambda}_{X_F}.
\end{equation}
Similarly, we consider  a multi-parameter family of ${\rm
\mathcal{K}}(n)$-modules on the space
$\mathbb{D}^n_{\lambda_1,\dots,\lambda_m;\mu}$ of multi-linear
differential operators: $~A: {\mathbb
F}_{\lambda_1}^n\otimes\cdots\otimes\mathbb{F}_{\lambda_m}^n\longrightarrow{\mathbb
F}_\mu^n$ with the natural $\mathcal{K}(n)$-action:
\begin{equation*}
{X_F}\cdot A={\mathbb L}_{X_F}^\mu\circ A-(-1)^{|A||F|}A\circ
{\mathbb L}_{X_F}^{\lambda_1,\dots,\lambda_m},\end{equation*}where
${\mathbb L}_{X_F}^{\lambda_1,\dots,\lambda_m}$ is defined by the
Leibnitz rule. We also consider the ${\rm \mathcal{K}}(n)$-module
$\Pi\left(\mathbb{D}^n_{\lambda_1,\dots,\lambda_m;\mu}\right)$
with the $\mathcal{K}(n)$-action ($\Pi$ is the change of parity
operator):
\begin{equation*}
{X_F}\cdot \Pi(A)=\Pi\left({\mathbb L}_{X_F}^\mu\circ
A-(-1)^{(|A|+1)|F|}A\circ {\mathbb
L}_{X_F}^{\lambda_1,\dots,\lambda_m}\right).\end{equation*} Since
$-\eta_i^2=\partial_x,$ and $\partial_i=\eta_i-\theta_i\eta_i^2,$
every differential operator $A\in\mathbb{D}^n_{\lambda,\mu}$ can
be expressed in the form
\begin{equation}
\label{diff1}
A(F\alpha_n^\lambda)=\sum_{\ell=(\ell_1,\dots,\ell_n)}a_\ell(x,\theta)
\eta_1^{\ell_1}\dots\eta_n^{\ell_n}(F)\alpha_n^\mu,
\end{equation}
where the coefficients $a_\ell(x,\theta)$ are arbitrary polynomial
functions.

The Lie superalgebra ${\mathcal K}(n-1)$ can be realized  as a
subalgebra of ${\mathcal K}(n)$:
\begin{equation*}
 {\mathcal K}(n-1)=\Big\{X_F\in{\mathcal
 K}(n)~|~\partial_{n}F=0\Big\}.
\end{equation*} Therefore,
$\mathbb{D}^n_{\lambda_1,\dots,\lambda_m;\mu}$ and
$\mathbb{F}^n_{\lambda}$ are ${\mathcal K}(n-1)$-modules. Note
also that, for any  $i$ in $\{1,2,\dots,n-1\}$, ${\mathcal
K}(n-1)$ is isomorphic to
\begin{equation*}
 {\mathcal K}(n-1)^i=\Big\{X_F\in{\mathcal
 K}(n)~|~\partial_{i}F=0\Big\}.
 \end{equation*}
 \begin{proposition}
\label{iso1} As a ${\mathcal K}(n-1)$-module, we have
\begin{align}\label{nizar}
\mathbb{D}^n_{\lambda,\mu;\nu}\simeq\widetilde{\mathbb{D}}^{n-1}_{\lambda,\mu;\nu}:=~~&
{\mathbb{D}}^{n-1}_{\lambda,\mu;\nu}\oplus
{\mathbb{D}}^{n-1}_{\lambda+\frac{1}{2},\mu+\frac{1}{2};\nu}\oplus
{\mathbb{D}}^{n-1}_{\lambda,\mu+\frac{1}{2};\nu+\frac{1}{2}}\oplus
{\mathbb{D}}^{n-1}_{\lambda+\frac{1}{2},\mu;\nu+\frac{1}{2}}\oplus\notag\\[6pt]
&\,
\Pi\left({\mathbb{D}}^{n-1}_{\lambda,\mu;\nu+\frac{1}{2}}\oplus
{\mathbb{D}}^{n-1}_{\lambda,\mu+\frac{1}{2};\nu}\oplus
{\mathbb{D}}^{n-1}_{\lambda+\frac{1}{2},\mu;\nu}\oplus
{\mathbb{D}}^{n-1}_{\lambda+\frac{1}{2},\mu+\frac{1}{2};\nu+\frac{1}{2}}
\right).
\end{align}
\end{proposition}
\begin{proofname}. For any $F\in\mathbb{R}[x,\theta]$, we write $$F=F_1+F_2\theta_n\quad \text{where}\quad\partial_nF_1=\partial_nF_2=0$$ and we prove that
$$\mathbb{L}^\lambda_{X_H}F=\mathbb{L}^\lambda_{X_H}F_1+(\mathbb{L}^{\lambda+\frac{1}{2}}_{X_H}F_2)\theta_n.$$
Thus, it is clear that the map
\begin{equation}\label{varphi}
\begin{array}{lcll} \varphi_\lambda:&\mathbb{F}^n_{\lambda} &\rightarrow&
\mathbb{F}^{n-1}_{\lambda}\oplus\Pi(\mathbb{F}^{n-1}_{\lambda+\frac{1}{2}})\\
&F\alpha_n^\lambda&\mapsto
&(F_1\alpha_{n-1}^\lambda,\,\Pi(F_2\alpha_{n-1}^{\lambda+\frac{1}{2}})),
\end{array}
\end{equation}
is ${\mathcal K}(n-1)$-isomorphism.  So, we get the natural
${\mathcal K}(n-1)$-isomorphism  from
$\mathbb{F}^n_{\lambda}\otimes\mathbb{F}^n_{\mu}$ to
\begin{equation*}
 \mathbb{F}^{n-1}_{\lambda}\otimes\mathbb{F}^{n-1}_{\mu}\oplus
\mathbb{F}^{n-1}_{\lambda}\otimes\Pi(\mathbb{F}^{n-1}_{\mu+\frac{1}{2}})
\oplus\Pi(\mathbb{F}^{n-1}_{\lambda+\frac{1}{2}})\otimes\mathbb{F}^{n-1}_{\mu}
\oplus\Pi(\mathbb{F}^{n-1}_{\lambda+\frac{1}{2}})\otimes\Pi(\mathbb{F}^{n-1}_{\mu+\frac{1}{2}})
\end{equation*} denoted $\psi_{\lambda,\mu}$.
Therefore, we deduce a ${\mathcal K}(n-1)$-isomorphism:
\begin{equation}\label{Psi}
\begin{array}{llll} \Psi_{\lambda,\mu,\nu}:&\widetilde{\mathbb{D}}^{n-1}_{\lambda,\mu;\nu}
 &\rightarrow&
{\mathbb{D}}^{n}_{\lambda,\mu;\nu}\\[2pt]
&A&\mapsto&\varphi^{-1}_\nu\circ A\circ\psi_{\lambda,\mu}.
\end{array}
\end{equation}Here, we identify
the ${\mathcal K}(n-1)$-modules via the following isomorphisms:
\begin{gather*}\begin{array}{llllllll}
\Pi\left({\mathbb{D}}^{n-1}_{\lambda,\mu;\nu'}\right)
&\rightarrow&
\mathrm{Hom_{diff}}\left(\mathbb{F}^{n-1}_\lambda\otimes\mathbb{F}^{n-1}_\mu,
\Pi(\mathbb{F}^{n-1}_{\nu'})\right), &\Pi(A)&\mapsto&\Pi\circ
A,\\[10pt] \Pi\left({\mathbb{D}}^{n-1}_{\lambda,\mu';\nu}\right)
&\rightarrow&
\mathrm{Hom_{diff}}\left(\mathbb{F}^{n-1}_\lambda\otimes
\Pi(\mathbb{F}^{n-1}_{\mu'}), \mathbb{F}^{n-1}_{\nu}\right),
&\Pi(A)&\mapsto& A\circ(1\otimes\Pi),\\[10pt]
\Pi\left({\mathbb{D}}^{n-1}_{\lambda',\mu;\nu}\right)
&\rightarrow&
\mathrm{Hom_{diff}}\left(\Pi(\mathbb{F}^{n-1}_{\lambda'})
\otimes\mathbb{F}^{n-1}_\mu, \mathbb{F}^{n-1}_{\nu}\right),
&\Pi(A)&\mapsto& A\circ(\Pi\otimes\sigma),\\[10pt]
\Pi\left({\mathbb{D}}^{n-1}_{\lambda',\mu';\nu'}\right)
&\rightarrow&
\mathrm{Hom_{diff}}\left(\Pi(\mathbb{F}^{n-1}_{\lambda'})\otimes\Pi(\mathbb{F}^{n-1}_{\mu'}),
\Pi(\mathbb{F}^{n-1}_{\nu'})\right), &\Pi(A)&\mapsto&\Pi\circ
A\circ(\Pi\otimes\sigma\circ\Pi),\\[10pt]
{\mathbb{D}}^{n-1}_{\lambda,\mu';\nu'} &\rightarrow&
\mathrm{Hom_{diff}}\left(\mathbb{F}^{n-1}_{\lambda}\otimes\Pi(\mathbb{F}^{n-1}_{\mu'}),
\Pi(\mathbb{F}^{n-1}_{\nu'})\right),&A&\mapsto&\Pi\circ
A\circ(1\otimes\Pi),\\[10pt]
{\mathbb{D}}^{n-1}_{\lambda',\mu';\nu} &\rightarrow&
\mathrm{Hom_{diff}}\left(\Pi(\mathbb{F}^{n-1}_{\lambda'})\otimes\Pi(\mathbb{F}^{n-1}_{\mu'}),
\mathbb{F}^{n-1}_{\nu}\right),&A&\mapsto&
A\circ(\Pi\otimes\sigma\circ\Pi),\\[10pt]
{\mathbb{D}}^{n-1}_{\lambda',\mu;\nu'} &\rightarrow&
\mathrm{Hom_{diff}}\left(\Pi(\mathbb{F}^{n-1}_{\lambda'})\otimes\mathbb{F}^{n-1}_{\mu},
\Pi(\mathbb{F}^{n-1}_{\nu'})\right),&A&\mapsto&\Pi\circ
A\circ(\Pi\otimes\sigma),
\end{array}
\end{gather*}
where
$\lambda'=\lambda+\frac{1}{2},~\mu'=\mu+\frac{1}{2},~\nu'=\nu+\frac{1}{2}$
and $\sigma(F)=(-1)^{|F|}F$. \end{proofname}\hfill$\Box$


\section{$\mathcal{K}(n)$-Invariant Binary Differential
Operators} In this section, we will classify all
$\mathcal{K}(n)$-invariant binary differential operators acting on
the spaces of weighted densities on $\mathbb{R}^{1|n}$ for
$n\geq2$. As a first step towards these classifications, we shall
need the list of binary $\mathcal{K}(1)$-invariant differential
operators acting on the spaces of weighted densities on
$\mathbb{R}^{1|1}$.
\subsection{$\mathcal{K}(1)$-invariant binary
differential operators}
In \cite{DYA}, Leites et al. classified all binary
$\mathcal{K}(1)$-invariant differential operators
$$
{\mathbb{F}}^1_\lambda\otimes{\mathbb{F}}^1_\mu\rightarrow{\mathbb{F}}^1_{\nu,\quad}F\alpha_1^\lambda\otimes
G\alpha_1^\mu\mapsto
\mathbb{T}_{\lambda,\mu,\nu}(F,G)\alpha_1^{\nu}.
$$
Recall that their list consists of (here
$\nu_k=\lambda+\mu+\frac{k}{2}$ for $k=0,1,2,3$)
\begin{equation} \label{k1}
\begin{array}{llllllllllll}
\mathbb{T}_{\lambda,\mu,\nu_0}(F,G)&=&
FG,\\[5pt]
\mathbb{T}_{0,0,\frac{1}{2}}^{a,b}(F,G)&=&a(-1)^{|F|}
F\eta_1(G)+b\eta_1(F)G,\quad a,\,b\in\mathbb{R},
\\[5pt]
\mathbb{T}_{\lambda,\mu,\nu_1}(F,G)&=&\mu\,\eta_1(F)G -\lambda
(-1)^{|F|}F\eta_1(G),\\[5pt]
\mathbb{T}_{\lambda,\mu,\nu_2}(F,G)&=&\mu F'G-{1\over
2}(-1)^{|F|}\eta_1(F)\eta_1(G)-\lambda
FG',\\[5pt]
\mathbb{T}_{0,\mu,\nu_3}(F,G)&=&S(F,G)
-2\mu\eta_1(F')G,\\[5pt]
\mathbb{T}_{\lambda,0,\nu_3}(F,G)&=&S(F,G)
-2\lambda(-1)^{|F|}F\eta_1(G',\\[5pt]
\mathbb{T}_{0,0,2}(F,G)&=&F'G'+(-1)^{|F|}\left(\eta_1(F')\eta_1(G)-\eta_1(F)\eta_1(G')\right),\\[5pt]
\mathbb{T}_{-\frac{3}{2},0,\frac{1}{2}}(F,G)&=&3FG''-(-1)^{|F|}M(F,G)
+2F'G',\\[5pt]
\mathbb{T}_{0,-\frac{3}{2},\frac{1}{2}}(F,G)&=&3F''G+(-1)^{|F|}M(G,F)
+2F'G',\\[5pt]
\mathbb{T}_{\lambda,-\lambda-1,\frac{1}{2}}(F,G)&=&\lambda(-1)^{|F|}
F\eta_1(G')+(\lambda+1)\eta_1(F')G+(\lambda+\frac{1}{2}) S(F,G),
\end{array}
\end{equation}
where
$$M(F,G)=2\eta_1(F)\eta_1(G')+\eta_1(F')\eta_1(G)\,\,\text{and}\,\,S(F,G)=\eta_1(F)G'+(-1)^{|F|}F'\eta_1(G).$$

Observe that the operation $\mathbb{T}_{\lambda,\mu,\nu_2}$ is
nothing but the well-known Poisson bracket on $\mathbb{R}^{1|1}$
 and the operation
$\mathbb{T}_{\lambda,\mu,\nu_1}$ is just the Buttin  bracket in
coordinates $\theta$ and $p:=\Pi(\a_1)$ with $x$ serving as
parameter (see, e.g. \cite{gmo, go, L, DYA}).

 \subsection{ $\mathcal{K}(n)$-invariant binary differential
operators for $n\geq2$}
Now, we describe the spaces of $\mathcal{K}(n)$-invariant binary
differential operators $\mathbb {F}^n_{\lambda}\otimes\mathbb
{F}^n_{\mu}\longrightarrow\mathbb {F}^n_{\nu}$ for $n\geq2$. We
prove that these spaces are nontrivial only if $\nu=\lambda+\mu$
or $\nu=\lambda+\mu+1$ and they are, in some sense,  spanned by
the following even operators defined on
$\mathbb{R}[x,\theta]\otimes\mathbb{R}[x,\theta]$:
\begin{equation}\label{invariants}
\begin{array}{llllllll} \mathfrak{a}(F,G)&=&FG,\\[3pt]
\mathfrak{b}(F,G)&=&\mu F'G-\lambda FG'-{1\over
2}(-1)^{|F|}\sum_{i=1}^n
\eta_i(F)\eta_i(G),\\[3pt]
\mathfrak{c}(F,G)&=&(-1)^{|F|}(
\eta_1(F)\eta_2(G)-\,\eta_2(F)\eta_1(G))+2\mu \eta_2(\eta_1(F))G),\\[3pt]
\mathfrak{d}(F,G)&=&(-1)^{|F|}(
\eta_1(F)\eta_2(G)-\,\eta_2(F)\eta_1(G))+2\lambda
F\eta_2\eta_1(G)),\\[3pt]
\mathfrak{e}(F,G)&=&(-1)^{|F|}(\lambda+{1\over
2})\left(\eta_1(F)\eta_2(G)- \eta_2(F)\eta_1(G)\right)+\lambda
F\eta_1\eta_2(G)+ (\lambda+1)\eta_1(\eta_2(F))G.
\end{array}
\end{equation}
More precisely, we have
\begin{thm}
\label{main} Let $n\geq2$ and $$\mathbb
{F}^n_{\lambda}\otimes\mathbb {F}^n_{\mu}\longrightarrow \mathbb
{F}^n_{\nu},\quad F\alpha_n^\lambda\otimes G \alpha_n^\mu\mapsto
\mathfrak{T}^n_{\lambda,\mu,\nu}(F,G)\alpha_n^\nu$$ be a
nontrivial $\mathcal{K}(n)$-invariant binary differential
operator. Then
$$\nu=\lambda+\mu\quad\text{or}\quad\nu=\lambda+\mu+1.$$
Moreover,

(a) If $\nu=\lambda+\mu$ then $\mathfrak{T}^n_{\lambda,\mu,\nu}=\alpha\mathfrak{a}.$

(b) If $\nu=\lambda+\mu+1$ then, for $n>2$ or $n=2$ but $\lambda\mu\nu\neq0$ we have $\mathfrak{T}^n_{\lambda,\mu,\nu}=\alpha\mathfrak{b}$ and if $n=2$ and $\lambda\mu\nu=0$  then $\mathfrak{T}^n_{\lambda,\mu,\nu}$  has the form $\alpha\mathfrak{b}+\beta\mathfrak{c}$, $\alpha\mathfrak{b}+\beta\mathfrak{d}$ or $\alpha\mathfrak{b}+\beta\mathfrak{e}$ in accordance with $\lambda=0$, $\mu=0$ or $\nu=0$. Here, $\alpha,\,\beta\in\mathbb{R}$ and $\mathfrak{a},\, \mathfrak{b},\, \mathfrak{c},\, \mathfrak{d},\, \mathfrak{e}$ are defined by (\ref{invariants}).
\end{thm}
\begin{proofname} (i) First assume that $n=2$.
The $\mathcal{K}(2)$-invariance of any element of
$\mathbb{D}^{2}_{\lambda,\mu;\nu}$ is equivalent to invariance
with respect just to the vector fields $X_F\in{\mathcal
 K}(2)$ such that $\partial_{1}\partial_{2}F=0$ that generate
 $\mathcal{K}(2)$. That is, an element  of
$\mathbb{D}^{2}_{\lambda,\mu;\nu}$ is $\mathcal{K}(2)$-invariant
if and only if it is invariant with respect just to the two
subalgebras $\mathcal{K}(1)$ and $\mathcal{K}(1)^1$. Obviously,
the $\mathcal{K}(1)$-invariant elements of
$\Pi(\mathbb{D}^{1}_{\lambda,\mu;\nu})$ can be deduced from those
given in (\ref{k1}) by using  the following
$\mathcal{K}(1)$-isomorphism
\begin{equation}
\label{od} \mathbb{D}^{1}_{\lambda,\mu;\nu}\rightarrow
\Pi(\mathbb{D}^{1}_{\lambda,\mu;\nu}),\quad~A\mapsto
\Pi(A\circ(\sigma\otimes\sigma))
\end{equation}  Now, by isomorphism (\ref{Psi}) we
exhibit the $\mathcal{K}(1)$-invariant elements of
$\mathbb{D}^{2}_{\lambda,\mu;\nu}$. Of course, these elements are
identically zero if
$2(\nu-\mu-\lambda)\neq-1,\,0,\,1,\,2,\,3,\,4,\,5,\,6.$ More
precisely, any $\mathcal{K}(1)$-invariant element $\mathfrak{T}$ of
$\mathbb{D}^{2}_{\lambda,\mu;\nu}$ can be expressed as follows
$$\begin{array}{llll}
\mathfrak{T}&=&\displaystyle\sum_{j,\ell,k=0,1}
\Omega_{\lambda,\mu,\nu}^{j,\ell,k}\Psi_{\lambda,\mu,\nu}\left(\Pi^{j+\ell+k}
\left(\mathbb{T}_{\lambda+{j\over 2},\mu+{\ell\over 2},\nu+{k\over
2}}\circ(\sigma^{j+\ell+k}\otimes\sigma^{j+\ell+k})\right)\right)+\\&~&\displaystyle\sum_{j,\ell,k=0,1}
\Omega_{\lambda,\mu,\nu,b}^{j,\ell,k,a}\Psi_{\lambda,\mu,\nu}\left(\Pi^{j+\ell+k}
\left(\mathbb{T}^{a,b}_{\lambda+{j\over 2},\mu+{\ell\over
2},\nu+{k\over
2}}\circ(\sigma^{j+\ell+k}\otimes\sigma^{j+\ell+k})\right)\right)\end{array}
$$
where $\mathbb{T}_{\lambda+{j\over 2},\mu+{\ell\over
2},\nu+{k\over 2}}, \mathbb{T}^{a,b}_{\lambda+{j\over
2},\mu+{\ell\over 2},\nu+{k\over 2}}$ are defined by (\ref{k1}).
The coefficients $\Omega_{\lambda,\mu,\nu}^{j,\ell,k}$ and
$\Omega_{\lambda,\mu,\nu,b}^{j,\ell,k,a}$ are, a priori, arbitrary
constants, but the invariance of $\mathfrak{T}$ with respect
$\mathcal{K}(1)^1$ imposes some supplementary conditions over
these coefficients and determines thus completely the space of
$\mathcal{K}(2)$-invariant elements of
$\mathbb{D}^{2}_{\lambda,\mu;\nu}$. By a direct computation, we
get:
$$\begin{array}{lllllllll}\Omega_{\lambda,\mu,\lambda+\mu}^{0,0,0}&=&\Omega_{\lambda,\mu,\lambda+\mu}^{0,1,1}
&=&\Omega_{\lambda,\mu,\lambda+\mu}^{1,0,1},&~&\\[6pt]\Omega_{\lambda,\mu,\lambda+\mu+1}^{0,0,0}&=&
2\Omega_{\lambda,\mu,\lambda+\mu+1}^{1,1,0}
&=&\Omega_{\lambda,\mu,\lambda+\mu+1}^{0,1,1}&=&
 \Omega_{\lambda,\mu,\lambda+\mu+1}^{1,0,1}.
\end{array}$$
All other coefficients vanish except for $\nu=\lambda+\mu+1$ with
$\lambda\mu\nu=0$, in which case we have also the following
non-zero coefficients:
\begin{equation*}
\begin{array}{llllllllll}
\Omega_{\lambda,0,\nu}^{0,0,1}&=&-\frac{1}{2}\Omega_{\lambda,0,\nu}^{0,1,0}
&=&
\frac{2\l+1}{2}\Omega_{\lambda,0,\nu}^{1,0,0}&=&\frac{1}{2}\Omega_{\lambda,0,\nu}^{1,1,1}&\hbox{
for } \lambda\neq-\frac{1}{2}
,\\[6pt]
\Omega_{\lambda,0,\nu}^{0,0,1}&=&-\frac{1}{2}\Omega_{\lambda,0,\nu}^{0,1,0}
&=&-
\Omega_{\lambda,0,\nu,0}^{1,0,0,1}&=&\frac{1}{2}\Omega_{\lambda,0,\nu}^{1,1,1}&\hbox{
for } \lambda=-\frac{1}{2}
,\\[6pt]
\Omega_{0,\mu,\nu}^{0,0,1}&=&-\frac{2\mu+1}{2}\Omega_{0,\mu,\nu}^{0,1,0}
&=&
\frac{1}{2}\Omega_{0,\mu,\nu}^{1,0,0}&=&\frac{1}{2}\Omega_{0,\mu,\nu}^{1,1,1}&\hbox{
for } \mu\neq-\frac{1}{2}
,\\[6pt]
\Omega_{0,\mu,\nu}^{0,0,1}&=&-\Omega_{0,\mu,\nu,1}^{0,1,0,0}
&=&\frac{1}{2}\Omega_{0,\mu,\nu}^{1,0,0}&=&\frac{1}{2}\Omega_{0,\mu,\nu}^{1,1,1}&\hbox{
for } \mu=-\frac{1}{2}
,\\[6pt]
\Omega_{\lambda,\mu,0}^{0,0,1}&=&\Omega_{\lambda,\mu,0}^{0,1,0}
&=&
\Omega_{\lambda,\mu,0}^{1,0,0}&=&(2\l+1)\Omega_{\lambda,\mu,0}^{1,1,1}&\hbox{
for } \lambda\neq-\frac{1}{2},\\[6pt]
\Omega_{\lambda,\mu,0}^{0,0,1}&=&\Omega_{\lambda,\mu,0}^{0,1,0}
&=&
\Omega_{\lambda,\mu,0}^{1,0,0}&=&-2\Omega_{\lambda,\mu,0,1}^{1,1,1,1}&\hbox{
for } \lambda=-\frac{1}{2}.
\end{array}
\end{equation*}
Thus, we easily check that Theorem \ref{main} is proved for $n=2$.

(ii) Now, we assume that $n\geq3$ and then we proceed by
recurrence over $n$. First note that the
$\mathcal{K}(n)$-invariance of any element of
$\mathbb{D}^{n}_{\lambda,\mu;\nu}$ is equivalent to invariance
with respect just to the fields $X_F\in{\mathcal
 K}(n)$ such that $\partial_{1}\cdots\partial_{n}F=0$ that generate
 $\mathcal{K}(n)$. That is, an element  of
$\mathbb{D}^{n}_{\lambda,\mu;\nu}$ is $\mathcal{K}(n)$-invariant
if and only if it is invariant with respect to the subalgebras
$\mathcal{K}(n-1)$ and $\mathcal{K}(n-1)^i,~i=1,\dots,n-1$. Thus,
as before, we prove that our result holds for $n=3.$ Assume that
it holds for $n\geq3$. Then, by recurrence assumption and
isomorphism (\ref{Psi}), we deduce that any nontrivial
$\mathcal{K}(n)$-invariant element $\mathfrak{T}$ of
$\mathbb{D}^{n+1}_{\lambda,\mu;\nu}$
 only can appear if
$2(\nu-\mu-\lambda)=-1,\,0,\,1,\,2,\,3,\,4,$ and it has the
general following form:
$$
\mathfrak{T}=\displaystyle\sum_{j,\ell,k=0,1}
\Delta_{\lambda,\mu,\nu}^{j,\ell,k}\Psi_{\lambda,\mu,\nu}\left(\Pi^{j+\ell+k}
\left(\mathfrak{T}^n_{\lambda+{j\over 2},\mu+{\ell\over
2},\nu+{k\over
2}}\circ(\sigma^{j+\ell+k}\otimes\sigma^{j+\ell+k})\right)\right).$$
As before, the coefficients $\Delta_{\lambda,\mu,\nu}^{j,\ell,k}$
are, a priori, arbitrary constants, but the invariance of
$\mathfrak{T}$ with respect $\mathcal{K}(n)^i$, $i=1,\dots,n$,
shows that
$$\begin{array}{lllllllllll}
\Delta_{\lambda,\mu,\lambda+\mu}^{0,0,0}&=&
\Delta_{\lambda,\mu,\lambda+\mu}^{0,1,1}&=&\Delta_{\lambda,\mu,\lambda+\mu}^{1,0,1},
&~&\\[6pt]
\Delta_{\lambda,\mu,\lambda+\mu+1}^{0,0,0}&=&2\Delta_{\lambda,\mu,\lambda+\mu+1}^{1,1,0}
&=&\Delta_{\lambda,\mu,\lambda+\mu+1}^{0,1,1}&=&\Delta_{\lambda,\mu,\lambda+\mu+1}^{1,0,1}
\end{array}
$$
and all other coefficients are identically zero. Therefore, we
easily check that $\mathfrak{T}$ is expressed as in Theorem
\ref{main}.\hfill $\Box$
\end{proofname}
\subsection{Poisson superalgebra of weighted densities}
For $n\geq2$, the even operation
\begin{equation}
\label{Poisson}
\mathfrak{T}^n_{\lambda,\mu,\lambda+\mu+1}(F,G)=\mu F'G-\lambda
FG'-{1\over 2}(-1)^{|F|}\sum_{i=1}^n \eta_i(F)\eta_i(G)
\end{equation}
defines a structure of Poisson Lie superalgebra on
$\mathbb{R}^{1|n}.$ Indeed, consider the continuous sum (direct
integral) of all spaces $\mathbb{F}^n_\lambda$:
\begin{equation*}
\mathbb{F}^n=\cup_{\lambda\in\mathbb{R}}\mathbb{F}_\lambda^n.
\end{equation*}
The collection of the operations (\ref{Poisson}) defines a
bilinear map $\mathfrak{T}^n_{1}:
\mathbb{F}^n\otimes\mathbb{F}^n\rightarrow\mathbb{F}^n$. The
following statement can be checked directly.
\begin{proposition}
The operation $\mathfrak{T}^n_{1}$ satisfies the Jacobi and
Leibniz identities, then it equips the space $\mathbb{F}^n$ with a Poisson
superalgebra structure.
\end{proposition}
Note that this Proposition is a simplest generalization of a
result by Gargoubi and Ovsienko for $n=1$ (see \cite{go}).
\section{Cohomology}
Let us first
recall some fundamental concepts from cohomology theory~(see, e.g.,
\cite{c}). Let $\frak{g}=\frak{g}_{\bar 0}\oplus \frak{g}_{\bar 1}$
be a Lie superalgebra acting on a superspace $V=V_{\bar 0}\oplus
V_{\bar 1}$ and let $\mathfrak{h}$ be a subalgebra of
$\mathfrak{g}$. (If $\frak{h}$ is omitted it assumed to be $\{0\}$.)
The space of $\frak h$-relative $n$-cochains of $\frak{g}$ with
values in $V$ is the $\frak{g}$-module
\begin{equation*}
C^n(\frak{g},\frak{h}; V ) := \mathrm{Hom}_{\frak
h}(\Lambda^n(\frak{g}/\frak{h});V).
\end{equation*}
The {\it coboundary operator} $ \delta_n: C^n(\frak{g},\frak{h}; V
)\longrightarrow C^{n+1}(\frak{g},\frak{h}; V )$ is a
$\frak{g}$-map satisfying $\delta_n\circ\delta_{n-1}=0$. The
kernel of $\delta_n$, denoted $Z^n(\mathfrak{g},\frak{h};V)$, is
the space of $\frak h$-relative $n$-{\it cocycles}, among them,
the elements in the range of $\delta_{n-1}$ are called $\frak
h$-relative $n$-{\it coboundaries}. We denote
$B^n(\mathfrak{g},\frak{h};V)$ the space of $n$-coboundaries.

By definition, the $n^{th}$ $\frak h$-relative  cohomolgy space is
the quotient space
\begin{equation*}
\mathrm{H}^n
(\mathfrak{g},\frak{h};V)=Z^n(\mathfrak{g},\frak{h};V)/B^n(\mathfrak{g},\frak{h};V).
\end{equation*}
We will only need the formula of $\delta_n$ (which will be simply
denoted $\delta$) in degrees 0 and 1: for $v \in
C^0(\frak{g},\,\frak{h}; V) =V^{\frak h}$,~ $\delta v(g) : =
(-1)^{p(g)p(v)}g\cdot v$, where
\begin{equation*}
V^{\frak h}=\{v\in V~\mid~h\cdot v=0\quad\text{ for all }
h\in\frak h\},
\end{equation*}
and  for  $ \Upsilon\in C^1(\frak{g}, \frak{h};V )$,
\begin{equation*}\delta(\Upsilon)(g,\,h):=
(-1)^{|g||\Upsilon|}g\cdot
\Upsilon(h)-(-1)^{|h|(|g|+|\Upsilon|)}h\cdot
\Upsilon(g)-\Upsilon([g,~h])\quad\text{for any}\quad g,h\in
\frak{g}.
\end{equation*}
\subsection{The space $\mathrm{H^1_{diff}}({\mathcal
K}(n);\mathbb{D}^n_{\lambda,\mu})$}
In this  subsection, we will compute the first differential
cohomology spaces  ${\mathrm H}^1_{\rm diff}(\mathcal{K}(n);
\mathbb{D}^n_{\lambda,\mu})$ for $n\geq3$. Our main result is the
following:
\begin{thm}
\label{main1}
  The space $\mathrm{H^1_{diff}}({\mathcal
K}(n);\mathbb{D}^n_{\lambda,\mu})$  has the following structure:
\begin{equation}
\mathrm{H^1_{diff}}({\mathcal
K}(n);\mathbb{D}^n_{\lambda,\mu})\simeq \left\{
\begin{array}{llllll}
\mathbb{R}&\text{if}\quad\left\{\begin{array}{llllll}n=3\,\text{ and }\,\mu-\lambda=0,\,\frac{1}{2},\,\frac{3}{2},\\n=4\,\text{ and }\, \mu-\lambda=0,\,1,\\n\geq5\,\text{ and }\mu-\lambda=0,\end{array}\right.\\[2pt]
0&\text{otherwise}.
\end{array}
\right.
\end{equation}
A base for the nontrivial $\mathrm{H^1_{diff}}({\mathcal K}(n);\mathbb{D}^n_{\lambda,\mu})$ is given by the cohomology
classes of the 1-cocycles:
\renewcommand{\arraystretch}{1.4}
$$
\begin{array}{lll}
\Upsilon^n_{\lambda,\lambda}(X_G)=G'
\\
\Upsilon^3_{\lambda,\lambda+\frac{1}{2}}(X_G)=
\left\{\begin{array}{ll}\eta_3\eta_2\eta_1(G)&\hfill\text{
 if }\lambda\neq-\frac{1}{2}\\
 \partial_{3}(G)
{\eta}_1{\eta}_2
-{\eta}_1{\eta}_2\left(\partial_{3}(G)\right)\zeta_4-
(-1)^{|G|}\theta_3M_{\eta_{3}(G)}{\eta}_{3} &\hfill\text{ if
}\lambda=-\frac{1}{2}\end{array}\right.\\
\Upsilon^3_{\lambda,\lambda+\frac{3}{2}}(X_G)=\left\{\begin{array}{ll}
\Xi_{G'}+2\lambda\eta_3\eta_2\eta_1(G')+\eta_3\eta_2\eta_1(G)\eta_1^2
~~~~~~~~~~~~~~~~~~&\hfill\text{
if }\lambda\neq-1\\[5pt]
\Xi_{G'}+\sum_{1\leq i<j\leq3}(-1)^{i+j}
\eta_{6-i-j}(G')\eta_j\eta_i &\hfill\text{ if }\lambda=-1
 \end{array}\right.\\
\Upsilon^4_{\lambda,\lambda+1}(X_G)= \left\{\begin{array}{ll}
Q_{G}+2\lambda\eta_4\eta_3\eta_2\eta_1(G)  &\hfill\text{
if }\lambda\neq-1\\[5pt]
A_G+\sum_{i=3}^4(-1)^{i}\theta_iM_{\eta_{i}\zeta_i(G)}{\eta}_i\partial_{7-i}+\\
+2\eta_1\eta_2\left(\partial_4\partial_3(G)\right)\zeta_4\zeta_3
+M_{\partial_4\partial_3(G)}\left(\zeta_3\zeta_4+\theta_3\theta_4\eta_4\eta_3\right)
&\hfill\text{ if }\lambda=-1,\end{array}\right. \\
\end{array}
$$
 where
 \begin{equation}
 \label{babba}
 \small{\begin{array}{l}
M_{G}=(-1)^{|G|}\sum_{i=1}^2(-1)^{i}{\eta}_{3-i}(G)\eta_i,
\qquad \Xi_{G}=(-1)^{|G|}\sum_{1\leq
i<j\leq3}(-1)^{i+j}\eta_j\eta_i(G)\eta_{6-i-j},\\
Q_{G}= (-1)^{|G|}\sum_{1\leq
i<j<k\leq4}(-1)^{i+j+k}\eta_k\eta_j\eta_i(G)\eta_{10-i-j-k},\\
A_G=(-1)^{|G|}\sum_{i=3}^4(-1)^{i}\Big({\eta}_1{\eta}_2\left(\partial_{i}\zeta_i(G)\right)\zeta_{7-i}
-\partial_{i}\zeta_i(G)
{\eta}_1{\eta}_2\Big)\partial_{7-i}\,\hbox{ with }\,\zeta_i=1-\theta_{7-i}\eta_{7-i}.
\end{array}}
 \end{equation}
\end{thm}

 The proof of Theorem \ref{main1} will be the subject of
subsection \ref{proof}. In fact, we need first the description of
$\mathrm{H^1_{diff}}({\mathcal
K}(n-1);\mathbb{D}^n_{\lambda,\mu})$ and
 $\mathrm{H^1_{diff}}({\mathcal
K}(n),{\mathcal K}(n-1)^i;\mathbb{D}^n_{\lambda,\mu})$.

\subsection{The space $\mathrm{H^1_{diff}}({\mathcal
K}(n-1);\mathbb{D}^n_{\lambda,\mu})$ }
The space $\mathrm{H^1_{diff}}({\mathcal
K}(n-1);\mathbb{D}^n_{\lambda,\mu})$ is closely related to
${\mathrm H}^1_{\rm diff}(\mathcal{K}(n-1);
\mathbb{D}^{n-1}_{\lambda,\mu})$. Therefore, for comparison and to
build upon, we first recall the description of ${\mathrm H}^1_{\rm
diff}(\mathcal{K}(2); \mathbb{D}^2_{\lambda,\mu})$. This space
was calculated in \cite{N}. The result is as follows:
\begin{equation}
\mathrm{H^1_{diff}}({\mathcal
K}(2);\mathbb{D}^2_{\lambda,\mu})\simeq\left\{
\begin{array}{llllll}
\mathbb{R}^2&\text{if}\quad \mu-\lambda=0,\,2,\\[2pt]
\mathbb{R}&\text{if}\quad \mu-\lambda=1,\\[2pt]
0&\text{otherwise}.
\end{array}
\right.
\end{equation}
The following 1-cocycles span the corresponding cohomology spaces:
\begin{equation}\label{maincocyc}
  \begin{array}{lllllllllll}
\Upsilon^2_{\lambda,\lambda}(X_G)&=&G'\\[2pt]
\widetilde{\Upsilon}^2_{\lambda,\lambda}(X_G)&=&\left\{\begin{array}{ll}\eta_1\eta_2(G)&\hfill\text{
 if }\lambda=0\\[2pt]2\,\lambda\,
\eta_1{\eta}_2\left(\theta_2\partial_{2}(G)\right) -(-1)^{|G|}
\sum_{i=1}^2{\eta}_i
\left(\theta_2\partial_{2}(G)\right){\eta}_{3-i} &\hfill\text{ if
}\lambda\neq0\end{array}\right. \\[2pt]
\Upsilon^2_{\lambda,\lambda+1}(X_G)&=&\left\{\begin{array}{ll}\eta_1\eta_2(G')
&\hfill\text{
if }\lambda\neq-{1\over2}\\[2pt]\eta_1\eta_2(G')+
 M_{G'}&\hfill\text{
if }\lambda=-{1\over2}
 \end{array}\right.\\[2pt]
 \Upsilon^2_{\lambda,\lambda+2}(X_G)&=&(2\lambda+1)\left(\frac{2\lambda}{3}G'''-H_{G''}\right)
-2 \eta_2\eta_1(G')\eta_2\eta_1
\\[2pt]
\widetilde{\Upsilon}^2_{\lambda,\lambda+2}(X_G)&=&\left\{\begin{array}{ll}
M_{G''}
+2\lambda\eta_2\eta_1(G'')-2\eta_2\eta_1(G')\partial_x&\hfill\text{if
}\lambda\neq-1\\[2pt]
\left(M_{G'}-\eta_2\eta_1(G')\right)\partial_x+M_{G''}-G''\eta_2\eta_1&\hfill\text{if
}\lambda=-1,\end{array}\right.
 \end{array}
\end{equation}
where, for $G\in\mathbb{R}[x,\theta],~M_{G}$ is as (\ref{babba})
and $H_G=(-1)^{|G|}\sum_{i=1}^2{\eta}_i(G)\eta_{i}$
 \begin{proposition}
\label{iso2} As a ${\mathcal K}(n-1)$-module, we have
\begin{equation}
\label{salem1}
\mathbb{D}^n_{\lambda,\mu}\simeq{\mathbb{D}}^{n-1}_{\lambda,\mu}\oplus
{\mathbb{D}}^{n-1}_{\lambda+\frac{1}{2},\mu+\frac{1}{2}}\oplus
\Pi\left({\mathbb{D}}^{n-1}_{\lambda,\mu+\frac{1}{2}}\oplus
{\mathbb{D}}^{n-1}_{\lambda+\frac{1}{2},\mu}\right).
\end{equation}
\end{proposition}
\begin{proofname}.
 By isomorphism (\ref{varphi}), we deduce a $\cK(n-1)$-isomorphism:
\begin{equation}\label{Phi}
\begin{array}{lll}&\Phi_{\lambda,\mu}:\mathbb{D}_{\lambda,\mu}^{n-1}\oplus
\mathbb{D}_{\lambda+\half,\mu+\half}^{n-1}\oplus
\Pi\left(\mathbb{D}_{\lambda,\mu+\half}^{n-1}\oplus
\mathbb{D}_{\lambda+\half,\mu}^{n-1}\right) \rightarrow
\mathbb{D}_{\lambda,\mu}^n\\[2pt]&A\mapsto\varphi_{\mu}^{-1}\circ
A\circ\varphi_{\lambda}.
\end{array}
\end{equation}Here, we identify
the $\cK(n-1)$-modules via the following isomorphisms:
\begin{gather*}\begin{array}{llllllll}
\Pi\left(\mathbb{D}^{n-1}_{\lambda,\mu+\frac{1}{2}}\right)&\rightarrow&
\mathrm{Hom_{diff}}\left(\mathbb{F}^{n-1}_\lambda,\Pi(\mathbb{F}^{n-1}_{\mu+\frac{1}{2}})\right)
\quad &\Pi(A)&\mapsto&\Pi\circ A,\\[10pt]
\Pi\left(\mathbb{D}^{n-1}_{\lambda+\frac{1}{2},\mu}\right)&\rightarrow&
\mathrm{Hom_{diff}}\left(\Pi(\mathbb{F}^{n-1}_{\lambda+\frac{1}{2}}),\mathbb{F}^{n-1}_{\mu}\right)
\quad &\Pi(A)&\mapsto& A\circ\Pi,\\[10pt]
\mathbb{D}^{n-1}_{\lambda+\frac{1}{2},\mu+\frac{1}{2}}&\rightarrow&
\mathrm{Hom_{diff}}\left(\Pi(\mathbb{F}^{n-1}_{\lambda+\frac{1}{2}}),\Pi(\mathbb{F}^{n-1}_{\mu+\frac{1}{2}})\right)
\quad &A&\mapsto&\Pi\circ A\circ\Pi.\\[10pt]
\end{array}
\end{gather*}
\end{proofname}.
\begin{corollary}
\label{main2} The space $\mathrm{H^1_{diff}}({\mathcal
K}(2);\mathbb{D}^3_{\lambda,\mu})$  has the following structure:
\begin{equation}
\mathrm{H^1_{diff}}({\mathcal
K}(2);\mathbb{D}^3_{\lambda,\mu})\simeq\left\{
\begin{array}{llllll}
\mathbb{R}^4&\text{if}\quad \mu-\lambda=0,\,2,\\[2pt]
\mathbb{R}^3&\text{if}\quad \mu-\lambda=\frac{1}{2},\,\frac{3}{2},\\[2pt]
\mathbb{R}^2&\text{if}\quad \mu-\lambda=-\frac{1}{2},\,1,\,\frac{5}{2},\\[2pt]
0&\text{otherwise}.
\end{array}
\right.
\end{equation}
The corresponding spaces $\mathrm{H^1_{diff}}({\mathcal
K}(2);\mathbb{D}^3_{\lambda,\lambda+{k\over2}})$ are spanned by
the cohomology classes of the 1-cocycles
$\Theta_{\lambda,\lambda+{k\over2}}^{3,j,\ell}$ and
$\widetilde{\Theta}_{\lambda,\lambda+{k\over2}}^{3,j,\ell},$
defined by
\begin{equation}
\label{Theta1} \Theta_{\lambda,\lambda+{k\over2}}^{3,j,\ell}(X_G)=
\Phi_{\lambda,\lambda+{k\over2}}\left(
\Pi^{j+\ell}\left(\sigma^{j+\ell}\circ\Upsilon^{2}_{\lambda+{j\over2},
\lambda+{k+\ell\over2}}(X_G)\right)\right)
\end{equation}
and
\begin{equation}
\label{Theta2}
\widetilde{\Theta}_{\lambda,\lambda+{k\over2}}^{3,j,\ell}(X_G)=
\Phi_{\lambda,\lambda+{k\over2}}\left(
\Pi^{j+\ell}\left(\sigma^{j+\ell}\circ\widetilde{\Upsilon}^{2}_{\lambda+{j\over2},
\lambda+{k+\ell\over2}}(X_G)\right)\right),
\end{equation}
where $j,\ell=0,1$, $k\in\{-1,\,\dots,\, 5\}$,
$\Upsilon^{2}_{\lambda,\mu},
\widetilde{\Upsilon}^{2}_{\lambda,\mu}$ are as in
(\ref{maincocyc}), $\Phi_{\lambda,\mu}$ is as in (\ref{Phi}).
Furthermore, the space $\mathrm{H^1_{diff}}({\mathcal
K}(2);\mathbb{D}^3_{\lambda,\lambda+{k\over2}})$ has the same
parity  as the integer $k$.
\end{corollary}
\begin{proofname}. First, it is easy to see that the map
$
\chi:\mathbb{D}^{n}_{\lambda,\mu}\rightarrow\Pi\left(\mathbb{D}^{n}_{\lambda,\mu}\right)
$ defined by $\chi(A)=\Pi(\sigma\circ A)$ satisfies
\begin{equation*}
\fL^{\lambda,\mu}_{X_G}\circ
\chi=(-1)^{|G|}\chi\circ\fL^{\lambda,\mu}_{X_G}\quad\text{ for all
} X_G\in\cK(n).
\end{equation*}
Thus, we deduce the structure of $\mathrm{H^1_{diff}}({\mathcal
K}(n);\Pi(\mathbb{D}^{n}_{\lambda,\mu}))$ from
$\mathrm{H^1_{diff}}({\mathcal
K}(n);\mathbb{D}^{n}_{\lambda,\mu})$. Indeed, to any $\Upsilon\in
Z^1_{\rm diff}({\mathcal K}(n);\mathbb{D}^{n}_{\lambda,\mu})$
corresponds $\chi\circ\Upsilon\in Z^1_{\rm diff}({\mathcal
K}(n);\Pi(\mathbb{D}^{n}_{\lambda,\mu}))$. Obviously, $\Upsilon$
is a couboundary if and only if  $\chi\circ\Upsilon$ is a
couboundary.

Second, according to Proposition \ref{iso2}, we obtain the
following isomorphism between cohomology spaces:
\begin{equation*}
\label{split}
\begin{array}{lll}
\HH^1_{\rm
diff}\left(\cK(n-1);\mathbb{D}^n_{\lambda,\mu}\right)\simeq&
\HH^1_{\rm
diff}\left(\cK(n-1);\mathbb{D}^{n-1}_{\lambda,\mu}\right)\oplus
\HH^1_{\rm
diff}\left(\cK(n-1);\mathbb{D}^{n-1}_{\lambda+\half,\mu+\half}\right)\oplus\\[2pt]
&\HH^1_{\rm
diff}\left(\cK(n-1);\Pi(\mathbb{D}^{n-1}_{\lambda+\half,\mu})\right)\oplus
\HH^1_{\rm
diff}\left(\cK(n-1);\Pi(\mathbb{D}^{n-1}_{\lambda,\mu+\half})\right).
\end{array}
\end{equation*}
Thus, we deduce the structure of $\mathrm{H^1_{diff}}({\mathcal
K}(2);\mathbb{D}^3_{\lambda,\mu})$.\end{proofname}
\hfill$\Box$

\subsection{The spaces $\mathrm{H^1_{diff}}({\mathcal
K}(n),{\mathcal K}(n-1)^i;\mathbb{D}^n_{\lambda,\mu})$ }
As a first step towards the proof of Theorem \ref{main1}, we shall
need to study the ${\mathcal K}(n-1)^i$-relative cohomology
$\mathrm{H^1_{diff}}({\mathcal K}(n),{\mathcal
K}(n-1)^i;\mathbb{D}^n_{\lambda,\mu})$. Hereafter all $\epsilon$'s
are constants and we will use the superscript $i$ when we consider
the superalgebra $\mathcal{K}(n)^i$ instead of $\mathcal{K}(n)$.

\medskip

Let $\mathfrak{g}=\mathfrak{h}\oplus\mathfrak{p}$ be a Lie
superalgebra, where $\mathfrak{h}$ is a subalgebra and
$\mathfrak{p}$ is a $\mathfrak{h}$-module such that
$[\mathfrak{p},\,\mathfrak{p}]=\mathfrak{h}$. Consider a 1-cocycle
$\Upsilon\in\mathrm{Z}^1(\mathfrak{g};\,V)$, where $V$ is a
$\mathfrak{g}$-module. The cocycle relation reads
$$
(-1)^{|g||\Upsilon|}g\cdot
\Upsilon(h)-(-1)^{|h|(|g|+|\Upsilon|)}h\cdot
\Upsilon(g)-\Upsilon([g,~h])=0\quad\text{ for any } g,h\in\mathfrak{g}.
$$
Denote $\Upsilon_{\mathfrak{h}}=\Upsilon|_{\mathfrak{h}}$ and
$\Upsilon_{\mathfrak{p}}=\Upsilon|_{\mathfrak{p}}$. Obviously, if
$\Upsilon_{\mathfrak{h}}=0$ then $\Upsilon$ is
$\mathfrak{h}$-invariant, therefore, the $\mathfrak{h}$-relative
cohomology space $\mathrm{H}^1(\mathfrak{g},\mathfrak{h};V)$ is
nothing but the space of cohomology classes of 1-cocycles
vanishing on $\mathfrak{h}$. In our situation,
$\mathfrak{g}={\mathcal K}(n),\,\mathfrak{h}={\mathcal
K}(n-1)^i,\,\mathfrak{p}=\Pi(\mathbb{F}^{n-1,i}_{-\frac{1}{2}})$
and $V=\mathbb{D}^n_{\lambda,\mu}$. Furthermore, in this case, the
1-cocycle relation yields the following equations:
\begin{align}
&\label{coc2}
(-1)^{|g||\Upsilon|}{X_{g}}\cdot\Upsilon_{\mathfrak{p}}(X_{h\theta_{i}})
-(-1)^{|\widetilde{h}|(|g|+|\Upsilon|)}
{X_{h\theta_{i}}} \cdot\Upsilon_{\mathfrak{h}}(X_{g})-\Upsilon_{\mathfrak{p}}([X_g,X_{h\theta_{i}}])=0,\\
&\label{coc3}(-1)^{|\widetilde{g}||\Upsilon|}{X_{g\theta_{i}}}\cdot\Upsilon_{\mathfrak{p}}(X_{h\theta_{i}})
-(-1)^{|\widetilde{h}|(|\widetilde{g}|+|\Upsilon|)}{X_{h\theta_{i}}}
\cdot\Upsilon_{\mathfrak{p}}(X_{g\theta_{i}})-\Upsilon_{\mathfrak{h}}([X_{g\theta_{i}},X_{h\theta_{i}}])=0,
\end{align}
where $g,
\,h\in\mathbb{R}[x,\theta_1,\dots,\breve{\theta}_i,\dots,\theta_n]$
and $|\widetilde{h}|=|h|+1$.
\begin{theorem}\label{relative}
For all $n\geq1$ and for all $i=1,\,\dots,\,n$, we have
\begin{equation}
\mathrm{H^1_{diff}}({\mathcal
K}(n),{\mathcal
K}(n-1)^i;\mathbb{D}^n_{\lambda,\mu})\simeq \left\{
\begin{array}{llllll}
\mathbb{R}&\text{if}\quad\left\{\begin{array}{llllll}n=2
\,\text{ and }\,\lambda=\mu\neq0,\\n=3\,\text{ and }\, (\lambda,\mu)=(-\frac{1}{2},0),\end{array}\right.\\[2pt]
0&\text{otherwise}.
\end{array}
\right.
\end{equation}
Moreover, for fixed $n=2$ or 3, non-zero relative cohomology
$\mathrm{H^1_{diff}}({\mathcal K}(n),{\mathcal
K}(n-1)^i;\mathbb{D}^n_{\lambda,\mu})$ are spanned by  classes of
some ${\mathcal K}(n-1)^i$-relative cocycles which are
cohomologous.
\end{theorem}
\begin{proofname}. For $n=1$, the result holds from \cite{bbbbk}, (Lemma 4.1).
For $n=2$, we deduce the result from \cite{N} (Proposition 4.2).
Moreover, the space $\mathrm{H^1_{diff}}({\mathcal K}(2),{\mathcal
K}(1);\mathbb{D}^2_{\lambda,\lambda})$, for $\lambda\neq0$, is
spanned by the cohomology class of the ${\mathcal K}(1)$-relative
1-cocycle $\widetilde{\Upsilon}_{\lambda,\lambda}$ defined by
(\ref{maincocyc}). Note that
$\widetilde{\Upsilon}_{\lambda,\lambda}|_{{\mathcal K}(1)^1}$ is a
coboundary, namely, $
2\delta\left(\theta_2\partial_{1}+\theta_1\partial_{2}\right)$.
Therefore, $\widetilde{\Upsilon}_{\lambda,\lambda}$ is
cohomologous to the ${\mathcal K}(1)^1$-relative 1-cocycle $\widetilde{\Upsilon}_{\lambda,\lambda}-2\delta\left(\theta_2\partial_{1}+\theta_1\partial_{2}\right)$
generating the space $\mathrm{H^1_{diff}}({\mathcal
K}(2),{\mathcal K}(1)^1;\mathbb{D}^2_{\lambda,\lambda})$.  Now, we
deduce the result for $n\geq3$ from the following Proposition.
\end{proofname}

\begin{prop}
\label{pro}~~
\begin{itemize}
  \item [1)] a)~For $(\lambda,\mu)\neq(-{1\over2},0),$ any element of
$Z_{\rm diff}^1(\cK(3);\mathbb{D}_{\lambda,\mu}^3)$ is a
coboundary over $\cK(3)$ if and only if at least one of its
restrictions to the subalgebras $\cK(2)^i$
 is a
coboundary.\\
  b)~ For $(\lambda,\mu)=(-{1\over2},0),$ there exists a unique, up to a scalar factor and
a coboundary, nontrivial 1-cocycle $\Upsilon^3_{-{1\over2},0}\in
Z_{\rm diff}^1(\cK(3);\mathbb{D}^3_{-{1\over2},0})$ such that its
restrictions to $\cK(2), ~\cK(2)^1$ and to $\cK(2)^2$ are
coboundaries. This 1-cocycle is odd and it is given by:
\begin{equation} \label{exp}
\Upsilon^3_{-{1\over2},0}(X_G)=
\partial_{3}(G)
{\eta}_1{\eta}_2
-{\eta}_1{\eta}_2\left(\partial_{3}(G)\right)\left(1-\theta_3\partial_{3}\right)-
(-1)^{|G|}\theta_3M_{\eta_{3}(G)}{\eta}_{3},
\end{equation}
where, for $G\in\mathbb{R}[x,\theta],~M_{G}$ is as (\ref{babba}).

\item [2)]~For $n>3,$ any element of
$Z_{\rm diff}^1(\cK(n);\mathbb{D}_{\lambda,\mu}^n)$ is a
coboundary over $\cK(n)$ if and only if at least one of its
restrictions to the subalgebras $\cK(n-1)^i$
is a coboundary.
\end{itemize}
\end{prop}
\begin{proofname}.
Let $\Upsilon\in Z_{\rm diff}^1(\cK(n);\mathbb{D}_{\lambda,\mu}^n)$
and assume that the restriction of
$\Upsilon$ to some $\mathcal{K}(n-1)^i$ is a coboundary, that
is, there exists $b\in \mathbb{D}_{\lambda,\mu}^n$ such that
$$
\Upsilon({X_F}) =\delta(b)(X_F)=(-1)^{|F||b|}X_F\cdot b
\quad\text{ for all}\quad
X_F\in\mathcal{K}(n-1)^i.
$$
By replacing $\Upsilon$ by
$\Upsilon-\delta{b}$, we can suppose that
$\Upsilon|_{\mathcal{K}(n-1)^i}=0$. Thus, the map $\Upsilon$ is
$\mathcal{K}(n-1)^i$-invariant and therefore the equation (\ref{coc3}) becomes:
\begin{equation}
 \label{sltr1}
 (-1)^{|\widetilde{g}||\Upsilon|}{X_{g\theta_{i}}}\cdot\Upsilon(X_{h\theta_{i}})
-(-1)^{|\widetilde{h}|(|\widetilde{g}|+|\Upsilon|)}{X_{h\theta_{i}}}
\cdot\Upsilon(X_{g\theta_{i}})=0.
\end{equation}
According to the isomorphism (\ref{varphi}), the map
$\Upsilon$ is decomposed into four components
\begin{equation}
\begin{array}{lllllllll}
&\Pi(\mathbb{F}^{n-1,i}_{-{1\over
2}})\otimes\mathbb{F}^{n-1,i}_\lambda
&\rightarrow&\mathbb{F}^{n-1,i}_\mu,\quad
&\Pi(\mathbb{F}^{n-1,i}_{-{1\over
2}})\otimes\Pi(\mathbb{F}^{n-1,i}_{\lambda+ {1\over
2}})&\rightarrow&\Pi(\mathbb{F}^{n-1,i}_{\mu+{1\over 2}}),\\
&\Pi(\mathbb{F}^{n-1,i}_{-{1\over
2}})\otimes\mathbb{F}^{n-1,i}_\lambda
&\rightarrow&\Pi(\mathbb{F}^{n-1,i}_{\mu+{1\over 2}}),\quad
&\Pi(\mathbb{F}^{n-1,i}_{-{1\over
2}})\otimes\Pi(\mathbb{F}^{n-1,i}_{\lambda+ {1\over
2}})&\rightarrow&\mathbb{F}^{n-1,i}_{\mu}.
\end{array}
\end{equation}
So, each of these bilinear maps
is $\mathcal{K}(n-1)^i$-invariant. Therefore, their expressions
are given by Theorem \ref{main}
 with the help of isomorphisms (\ref{Psi}) and (\ref{od}). More
precisely, using equation (\ref{sltr1}), we get up to a scalar
factor:

$\bullet$ For $n\geq3$ with $(\lambda,\mu)\neq(-{1\over 2},0)$ if
$n=3,$
\begin{equation}
\label{vanishing} \Upsilon=\left\{\begin{array}{lllll}
\delta\left(\theta_i\right)&\text{ if
}\quad\mu=\lambda-{1\over 2},\\[2pt] \delta
\left(1-2\theta_i\partial_{{i}}\right)&\text{
if }\quad \mu=\lambda,\\[2pt]
\delta\left(\partial_{{i}}\right)&\text{ if
}\quad\mu=\lambda+{1\over 2},\\[2pt]
 0 &\text{ otherwise. }
\end{array}\right.
\end{equation}

$\bullet$ For $n=3$ and $\mu=\lambda+{1\over 2}=0,$
\begin{equation}
\label{m} {\Upsilon}=\epsilon_1\,\Upsilon^{3, i}_{-{1\over
2},0}+\epsilon_2\,\delta \left(\partial_{i}\right),
\end{equation}
where $\Upsilon^{3, i}_{-{1\over 2},0}$ is the 1-cocycle on
$\cK(3)$ with coefficients in $\mathbb{D}^3_{-{1\over 2},0}$
defined by
\begin{equation} \label{i}
\begin{array}{lll}
\Upsilon^{3, i}_{-{1\over2},0}(X_G)&=&
\partial_{i}(G) {\eta}_\ell{\eta}_k
-{\eta}_\ell{\eta}_k\left(\partial_{i}(G)\right)
\left(1-\theta_i\partial_{i}\right) +
\\&&
~~+\theta_i\left({\eta}_\ell\eta_{i}(G){\eta}_k
-{\eta}_k\eta_{i}(G){\eta}_\ell \right){\eta}_{i}
\end{array}
\end{equation}
with $\ell,k\neq i$ and $\ell<k.$ Obviously, $\Upsilon^{3,
3}_{-{1\over2},0}=\Upsilon^{3}_{-{1\over 2},0}$ with
$\Upsilon^{3}_{-{1\over 2},0}$ is as in (\ref{exp}), and a direct
computation shows that for $j=1,2:$
\begin{equation}
\label{ni} \Upsilon^{3}_{-{1\over 2},0}+(-1)^j \Upsilon^{3,
j}_{-{1\over 2},0}= 2(-1)^j\delta\left((\theta_{3}\eta_{j}+
\theta_{j}\eta_{3})\eta_{3-j} \right).
\end{equation}
Thus, up to a scalar factor and a coboundary,
${\Upsilon}=\Upsilon^{3}_{-{1\over 2},0}.$ Therefore, in order to
complete the proof of Proposition \ref{pro}, we have to study the
cohomology class of the 1-cocycle $\Upsilon^{3}_{-{1\over 2},0}$
in $\HH^1_{\rm diff} (\cK(3) , \mathbb{D}^3_{-{1\over 2},0}).$
\begin{lem}
\label{nontrivial} The 1-cocycle $\Upsilon^{3}_{-{1\over 2},0}$
defines a nontrivial cohomology class over $\cK(3)$. Its
restrictions to $\cK(2),~\cK(2)^1$ and to $\cK(2)^2$ are
coboundaries.
\end{lem}
\begin{proofname}. It follows from equation (\ref{ni}) that the
restriction of $\Upsilon^{3}_{-{1\over 2},0}$ to $\cK(2)$ vanishes
and  to $\cK(2)^1$ and to $\cK(2)^2$ are coboundaries. Now, assume
that there exists an odd operator $A\in\mathbb{D}^3_{-{1\over
2},0}$ such that $\Upsilon^{3}_{-{1\over 2},0}$ is equal to
$\delta A.$ By isomorphism (\ref{Phi}), the operator $A$ can be
expressed as $ A=\Phi_{-{1\over 2},0}( A_{1},A_{2},\Pi(
A_{3}),\Pi( A_{4})), $ where $ A_{1}\in\mathbb{D}^2_{-{1\over
2},0}, A_{2}\in\mathbb{D}^2_{0,{1\over 2}},
A_{3}\in\mathbb{D}^2_{-{1\over 2},{1\over 2}}$ and $
A_{4}\in\mathbb{D}^2_{0,0}$. Thus, since the map
\begin{equation*}
\mathbb{D}^{2}_{\lambda,\mu}\rightarrow
\Pi(\mathbb{D}^{2}_{\lambda,\mu}),~B\mapsto \Pi(B\circ\sigma)
\end{equation*}
is a $\cK(2)$-isomorphism, the condition $\Upsilon^{3}_{-{1\over
2},0|\cK(2)}=0$ tell us that $A_{1},A_{2}, A_{3}\circ\sigma$ and
$A_{4}\circ\sigma$ are $\cK(2)$-invariant linear maps. Therefore,
up to a scalar factor, each of $A_{1},A_{2}, A_{3}\circ\sigma$ and
$A_{4}\circ\sigma$ is the identity  map \cite{L}: $
\mathbb{F}^2_\lambda\rightarrow\mathbb{F}^2_\lambda,\quad
F\alpha_{2}^\lambda\mapsto F\alpha_{2}^\lambda. $ Thus, we obtain
\begin{equation*}
A(F\alpha_2^{-{1\over 2}})=\epsilon\partial_3(F).
\end{equation*}
Finally, it is  easy to check that the equation
$\Upsilon^{3}_{-{1\over 2},0}=\delta(A)$ has no solutions
contradicting our assumption.
 Lemma \ref{nontrivial} is proved.
\end{proofname}
Thus we have completed the proof of Proposition \ref{pro}.
\hfill$\Box$
\end{proofname}
\begin{cor}\label{cor1}
\label{sd} Up to a coboundary, any 1-cocycle $\Upsilon\in
Z^1_\mathrm{diff}(\cK(3);\mathbb{D}^3_{\lambda,\mu})$ has the
following general form:
\begin{equation}\label{coef}
\Upsilon(X_{F})=\sum
a_{\ell_1\ell_2\ell_3 m_1m_2
m_3}\eta_1^{\ell_1}\eta_2^{\ell_2}\eta_3^{\ell_3}(F)\eta_1^{k_1}\eta_2^{k_2}\eta_3^{k_3},
\end{equation}
where the coefficients $a_{\ell_1\ell_2\ell_3 k_1k_2 k_3}$ are
functions of $\theta_i$, not depending on $x$.
\end{cor}

\begin{proofname}. By (\ref{diff1}), we can see that the operator $\Upsilon$ has
the form (\ref{coef}) where, a priori, the coefficients
$a_{\ell_1\ell_2\ell_3 k_1k_2 k_3}$ are some functions of $x$ and
$\theta_i$, but we shall now prove that
$\partial_x\,a_{\ell_1\ell_2\ell_3 k_1k_2 k_3}=0$. To do that, we
shall simply show that $X_1\cdot\Upsilon=0$.

We have
\begin{equation}\label{part}
(X_1\cdot\Upsilon)(X_{F}):= X_1\cdot\Upsilon(X_{F})-
\Upsilon([X_1,X_{F}])\quad\hbox{for all}~~ F\in
{\mathbb{R}}[x,\,\theta].
\end{equation}
But, from Proposition \ref{pro}
and Corollary
\ref{main2}, it follows that, up to a coboundary,
$\Upsilon(X_1)=0$, and therefore the equation (\ref{part}) becomes
\begin{equation}
\begin{array}{lll}\label{part1}
(X_1\cdot\Upsilon)(X_{F})&= X_1\cdot\Upsilon(X_{F})-
(-1)^{|F||\Upsilon|}{X_F}\cdot\Upsilon(X_1)- \Upsilon([X_1,X_F]).
\end{array}
\end{equation}
The right-hand side of (\ref{part1}) vanishes because $\Upsilon$
is a 1-cocycle. Thus, $X_1\cdot\Upsilon=0$.
\hfill$\Box$\end{proofname}

\medskip

The following lemma gives a description of all coboundaries over
$\mathcal{K}(2)^i$ vanishing on the subalgebra $\cK(1)^{m_i},$
where $m_i\in\{1,\,2,\,3\}\setminus\{i\}$. This
description will be useful in the proof of Theorem \ref{main1}.
\begin{lem} (see \cite{N})
\label{co}
Any coboundary $\mathcal{B}^{i,m_i}_{\lambda,\mu}\in
B^1_{\mathrm{diff}}(\cK(2)^i;\mathbb{D}^{2,i}_{\lambda,\mu})$
vanishing on $\cK(1)^{m_i}$ is, up to a scalar factor, given by
\begin{equation}
\label{coboundary}
 \mathcal{B}^{i,m_i}_{\lambda,\mu}=\left\{
 \begin{array}{llll}
 \delta  \left(\epsilon_1\,\partial_{m_i}+
\epsilon_2\eta_{6-i-m_i}( \theta_{m_i}\eta_{m_i}-1)\right)&\text{
if }
\quad(\lambda,\mu)=(0,\half)\\
\delta(\epsilon_1\,\partial_{m_i}+\epsilon_2\theta_{m_i}\eta_{6-i-m_i}\eta_{m_i})&\hbox{
if }\quad(\lambda,\mu)=(-\half,0)\\
\delta\left(\epsilon_1\,\theta_{m_i}\eta_{6-i-m_i}+
\epsilon_2\,\theta_{m_i}\eta_{{m_i}}\right)&\text{ if
}\quad\lambda=\mu=0\\
\delta\left(\partial_{m_i}\eta_{{6-i-m_i}}\right)&\text{
if }\quad(\lambda,\mu)=(-\half,\half)\\
\delta\left(\theta_{m_i}\eta_{{m_i}}\right)&\text{ if }\quad
\lambda=\mu\neq0\\
\delta \left(\partial_{m_i}\right) &\text{ if
}\quad\mu=\lambda+\half\hbox{ and }\lambda\neq0,-\half\\
\delta \left(\theta_{m_i}\right) &\text{ if
}\quad\mu=\lambda-\half\\
 0 &\hbox{ otherwise. }
\end{array}
\right.
\end{equation}
\end{lem}
\subsection{Proof of Theorem \ref{main1}}
\label{proof} (i) According to Proposition \ref{pro}, the
restriction of any nontrivial differential $1$-cocycle $\Upsilon$
of $\cK(3)$ with coefficients in $\mathbb{D}^3_{\lambda,\mu}$ to
$\cK(2)^i,$ for $i$=1, 2, 3, is a nontrivial 1-cocycle except for
\begin{equation}
\label{form} \Upsilon=\epsilon\Upsilon^{3}_{-{1\over 2},0}+\delta
A,
\end{equation}
 where $\Upsilon^{3}_{-{1\over 2},0}$ is as (\ref{exp}),
 $\epsilon\neq0$
and $A\in\mathbb{D}^3_{-{1\over 2},0}$. So, if
$2(\mu-\lambda)\neq-1,\,0,\,1,\,2,\,3,\,4,\,5$, then, by Corollary
\ref{main2}, the corresponding cohomology spaces ${\mathrm
H}^1_{\mathrm{diff}}(\cK(3);\mathbb{D}^3_{\l,\mu})$ vanish.

For $2(\mu-\lambda)=-1,\,0,\,1,\,2,\,3,\,4,\,5$, let $\Upsilon$ be
a 1-cocycle from $\cK(3)$ to $\mathbb{D}^3_{\l,\mu}$. The map
$\Upsilon_{|\cK(2)^i}$ is a 1-cocycle of $\cK(2)^i$. Therefore,
using Corollary \ref{main2} together with Lemma \ref{co} and
Theorem \ref{relative} with the help of isomorphism (\ref{Phi}),
we deduce that, up to a coboundary, the non-zero restrictions of
the cocycle $\Upsilon$ on $\cK(2)^i$ 
can be expressed as (here $\tau=\mu-\lambda$):

\noindent For $\tau=-\frac{1}{2},\, \frac{3}{2},\, 2$,
{\small\begin{equation*}
{\Upsilon}|_{\cK(2)^i}=\left\{
\begin{array}{llll}
a(-1)^{i}\widetilde{\Theta}_{\lambda,\mu}^{i,0,1}+\frac{i(3-i)}{2}b
\left(\Gamma_{\lambda,\mu}^{i,3,0,0}+\Gamma_{\lambda,\mu}^{i,3,1,1}\right)
&\text{if } \tau=-\frac{1}{2},\,\lambda=0\\[5pt]
a\left(\Lambda^i_{2,3}\delta(\widetilde{A}_{01})-(-1)^{i}\widetilde{\Theta}_{\lambda,\mu}^{i,0,1}\right)
+b({1\over2}\Lambda^i_{2,3}-\Lambda^i_{1,3})
\left(\Gamma_{\lambda,\mu}^{i,3,0,0}+\Gamma_{\lambda,\mu}^{i,3,1,1}\right)
&\text{if } \tau=-\frac{1}{2},\,\lambda\neq0\\[5pt]
a(-1)^{i}(\widetilde{\Theta}_{\lambda,\mu}^{i,0,1}-
\Theta_{\lambda,\mu}^{i,1,0})
&\text{if } \tau=\frac{3}{2}\\[5pt]
a\left((2\lambda+2)\Theta_{\lambda,\mu}^{i,0,0}+(2\lambda+3)\Theta_{\lambda,\mu}^{i,1,1}\right)
&\text{if }\tau=2\\[10pt]
\end{array}
\right.\end{equation*} }

\noindent For $\tau=0$, {\small\begin{equation*}
{\Upsilon}|_{\cK(2)^i}=\left\{
\begin{array}{llll}
a\left(\Theta_{\lambda,\mu}^{i,0,0}+\Theta_{\lambda,\mu}^{i,1,1}\right)+
b_i\left(\widetilde{\Theta}_{\lambda,\mu}^{i,0,0}+
\widetilde{\Theta}_{\lambda,\mu}^{i,1,1}\right)+\\
\Lambda^i_{1,3}\left(
2b_{1}\left(\Gamma_{\lambda,\mu}^{2,3,1,0}-\Gamma_{\lambda,\mu}^{2,3,0,1}\right)-
t\left(\Gamma_{\lambda,\mu}^{2,3,0,0}+\Gamma_{\lambda,\mu}^{2,3,1,1}\right)\right)+\\
\Lambda^i_{2,3} \left(b_{3}
\delta(\widetilde{A})+b_{2}(\Gamma_{\lambda,\mu}^{1,3,0,1}
-
\Gamma_{\lambda,\mu}^{1,3,1,0}) +{1\over
2}t\left(\Gamma_{\lambda,\mu}^{1,3,0,0}+\Gamma_{\lambda,\mu}^{1,3,1,1}\right)\right)&\text{
if } \mu\neq0,-{1\over2}\\[10pt]
a\left(\Theta_{\lambda,\mu}^{i,0,0}+\Theta_{\lambda,\mu}^{i,1,1}\right)
+b_i\widetilde{\Theta}_{\lambda,\mu}^{i,1,1} -2\Lambda^i_{1,3}
\Big(b_2\Gamma_{\lambda,\mu,1,0}^{2,3,0,0}+\\
b_3\Gamma_{\lambda,\mu,0,1}^{2,3,0,1}+
b_1(\Gamma_{\lambda,\mu,1,0}^{2,3,0,1}-
\Gamma_{\lambda,\mu}^{2,3,1,0})
+{1\over2}t(\Gamma_{\lambda,\mu}^{2,3,1,1}+\Gamma_{\lambda,\mu,0,1}^{2,3,0,0})\Big)+\\
\Lambda^i_{2,3} \Big(b_{3}( \delta(\widetilde{A}_{11})-
\Gamma_{\lambda,\mu,0,1}^{1,3,0,1})+
b_2(\Gamma_{\lambda,\mu,1,0}^{1,3,0,1}-
\Gamma_{\lambda,\mu}^{1,3,1,0})
+\\b_{1}\Gamma_{\lambda,\mu,1,0}^{1,3,0,0}+
{1\over2}t(\Gamma_{\lambda,\mu}^{1,3,1,1}+\Gamma_{\lambda,\mu,0,1}^{1,3,0,0})
\Big) &\text{
if } \mu=0\\[10pt]
a\left(\Theta_{\lambda,\mu}^{i,0,0}+\Theta_{\lambda,\mu}^{i,1,1}\right)
+b_i\widetilde{\Theta}_{\lambda,\mu}^{i,0,0} -2\Lambda^i_{1,3}
\Big(b_2\Gamma_{\lambda,\mu,1,0}^{2,3,1,1}+\\
b_3\Gamma_{\lambda,\mu,0,1}^{2,3,0,1}+
b_1(\Gamma_{\lambda,\mu,1,0}^{2,3,0,1}-
\Gamma_{\lambda,\mu}^{2,3,1,0})
+{1\over2}t(\Gamma_{\lambda,\mu}^{2,3,0,0}+\Gamma_{\lambda,\mu,0,1}^{2,3,1,1})\Big)+\\
\Lambda^i_{2,3} \Big(b_{3}( \delta(\widetilde{A}_{00})-
\Gamma_{\lambda,\mu,0,1}^{1,3,0,1})+
b_2(\Gamma_{\lambda,\mu,1,0}^{1,3,0,1}-
\Gamma_{\lambda,\mu}^{1,3,1,0})
+\\b_{1}\Gamma_{\lambda,\mu,1,0}^{1,3,1,1}+{1\over2}t(\Gamma_{\lambda,\mu}^{1,3,0,0}
+\Gamma_{\lambda,\mu,0,1}^{1,3,1,1}) \Big)
&\text{if } \mu=-{1\over2}\\[10pt]
\end{array}
\right. \end{equation*}}

\noindent For $\tau=\frac{1}{2}$, {\small\begin{equation*}
{\Upsilon}|_{\cK(2)^i}=\left\{
\begin{array}{llll}
b(-1)^{i}(\Theta_{\lambda,\mu}^{i,0,1}-\widetilde{\Theta}_{\lambda,\mu}^{i,1,0})+
({1\over2}\Lambda^i_{2,3}-\Lambda^i_{1,3})\times \\
~\times\left(a
(\Gamma_{\lambda,\mu,1,0}^{i,3,0,0}+\Gamma_{\lambda,\mu,1,0}^{i,3,1,1})
+(-1)^it(\Gamma_{\lambda,\mu,1,0}^{i,3,1,0}-\Gamma_{\lambda,\mu}^{i,3,0,1})\right)
&\text{if } \mu=0\\[10pt]
(-1)^{i}(a\Theta_{\lambda,\mu}^{i,0,1}+b\widetilde{\Theta}_{\lambda,\mu}^{i,1,0})+
t({1\over2}\Lambda^i_{2,3}-\Lambda^i_{1,3})\times \\~\times (
\Gamma_{\lambda,\mu,1,0}^{i,3,0,0}+\Gamma_{\lambda,\mu}^{i,3,1,1})-
b\Lambda^i_{2,3}\delta(\widetilde{A}_{10})
&\text{if } \mu=\frac{1}{2}\\[10pt]
(-1)^{i}(a\Theta_{\lambda,\mu}^{i,0,1}+b\widetilde{\Theta}_{\lambda,\mu}^{i,1,0})+
t({1\over2}\Lambda^i_{2,3}-\Lambda^i_{1,3})\times \\~\times (
\Gamma_{\lambda,\mu}^{i,3,0,0}+\Gamma_{\lambda,\mu,1,0}^{i,3,1,1})-
b\Lambda^i_{2,3}\delta(\widetilde{A}_{10})
&\text{if } \mu=-\frac{1}{2}\\[10pt]
(-1)^{i}(a\Theta_{\lambda,\mu}^{i,0,1}+b\widetilde{\Theta}_{\lambda,\mu}^{i,1,0})+
t({1\over2}\Lambda^i_{2,3}-\Lambda^i_{1,3})\times \\~\times (
\Gamma_{\lambda,\mu}^{i,3,0,0}+\Gamma_{\lambda,\mu}^{i,3,1,1})-
b\Lambda^i_{2,3}\delta(\widetilde{A}_{10})
&\text{if } \mu\neq\pm\frac{1}{2},0,\\[10pt]
\end{array}
\right. \end{equation*}}

\noindent where (recall that $\mathcal{B}^{i,m_i}_{\lambda, \mu}$
depend on $\epsilon_1,\epsilon_2$)
{\small \begin{equation*}
\begin{array}{lll}
 \Gamma_{\lambda,\mu,\epsilon_1,\epsilon_2}^{i,m_i,j,\ell}(X_G)=
 \Phi^i_{\lambda,\mu}\left(
\Pi^{j+\ell}\left(\sigma^{j+\ell}\circ\mathcal{B}^{i,m_i}_{\lambda+{j\over2},
\mu+{\ell\over2}}(X_G)\right)\right),\\
\widetilde{A}_{j\ell}=\Phi_{\lambda,\mu}\left(
\Pi^{j+\ell}\left(\sigma^{j+\ell}\circ A_{\lambda+{j\over2},
\mu+{\ell\over2}}\right)\right)\quad\hbox{with}\quad
A_{\lambda+{j\over2}, \mu+{\ell\over2}}=\theta_2\partial_{1}+
\theta_1\partial_{2}\in\mathbb{D}^{2}_{\l+{j\over2},\mu+{\ell\over2}},\\
\widetilde{A}=\widetilde{A}_{11}+\widetilde{A}_{00},\quad
\Lambda^i_{r,s}=(i-r)(i-s),
\end{array}
\end{equation*}}
\noindent $\Theta_{\lambda,\mu}^{i,j,\ell}$ and
$\widetilde{\Theta}_{\lambda,\mu}^{i,j,\ell}$ are defined by
(\ref{Theta1})--(\ref{Theta2}) and the coefficients $a,~b,~b_i$
and $t$ are constants. So, by Proposition \ref{pro}, $\HH^1_{\rm
diff} (\cK(3) ; \mathbb{D}^3_{\lambda,\mu})=0$ for $\mu-\lambda=1,
\,{5\over2}$. Now, by Corollary \ref{cor1}, we can write {\small
\begin{equation*}
\Upsilon(X_{h\theta_1\theta_2\theta_3})=
\sum_{m,k,~\varepsilon=(\varepsilon_1,\varepsilon_2,\varepsilon_3)}\left(a_{0,m,k,\varepsilon}+
\sum_{j=1}^3\sum_{1\leq i_1<\cdots <i_j\leq3}a_{i_1\cdots
i_j,m,k,\varepsilon}\theta_{i_1}\cdots\theta_{i_j}\right)h^{(k)}\partial_x^{m}
\partial_{1}^{\varepsilon_1}\partial_{2}^{\varepsilon_2}\partial_{3}^{\varepsilon_3}
\end{equation*}}

\noindent with $\varepsilon_i=0,\,1$. For each case, we solve the
equations (\ref{coc2}) and (\ref{coc3}) for $a,\,b,\,b_i,\,t,\,
a_{0,m,k,\varepsilon},$ $a_{i_1\cdots i_j,m,k,\varepsilon}$. We
obtain
\begin{itemize}
\item[1)] For $2(\mu-\lambda)=-1,\,4,$ the coefficient $a$
vanishes; so, by Proposition \ref{pro}, $\Upsilon$ is a
coboundary. Hence $\HH^1_{\rm diff} (\cK(3) ;
\mathbb{D}^3_{\lambda,\mu})=0$.
\item[2)] For $\mu=\lambda,$ the coefficients $b_i$ vanish and, up to a coboundary, $\Upsilon$ is a
multiple of $\Upsilon^3_{\lambda,\lambda},$ see Theorem
\ref{main1}. Hence \rm{dim}$\HH^1_{\rm diff} (\cK(3);
\mathbb{D}^3_{\lambda,\lambda})=1$.
\item[3)] For $2(\mu-\lambda)=1,$ the coefficient $b$ vanishes  and, up to a coboundary,
$\Upsilon$ is a multiple of
$\Upsilon^3_{\lambda,\lambda+{1\over2}},$ see Theorem \ref{main1}.
Hence \rm{dim}$\HH^1_{\rm diff} (\cK(3) ;
\mathbb{D}^3_{\lambda,\lambda+{1\over2}})=1$.
\item[4)] For $2(\mu-\lambda)=3$, $\Upsilon$ is a
multiple of $\Upsilon^3_{\lambda,\lambda+{3\over2}}$. Hence
$\hbox{\rm{dim}}\HH^1_{\rm diff} (\cK(3) ;
\mathbb{D}^3_{\lambda,\lambda+{3\over2}})=1$.
\end{itemize}

(ii) Note that, by Proposition \ref{pro}, the
restriction of any nontrivial differential $1$-cocycle $\Upsilon$
of $\cK(4)$ with coefficients in $\mathbb{D}^4_{\lambda,\mu}$ to
$\cK(3)^i,$ for $i=1,\ldots,4$, is a nontrivial 1-cocycle.
Furthermore, using arguments similar to those of the proof of
Corollary \ref{main2} together with the above result,
we deduce that  $\HH^1_{\rm diff} (\cK(3)^i ;
\mathbb{D}^4_{\lambda,\mu})=0$  if $2(\mu-\lambda)\neq-1,\,0,\,1,\,2,\,3,\,4.$ Then, we consider only the cases where $2(\mu-\lambda)=-1,\,0,\,1,\,2,\,3,\,4$ and, as before, we get the result for $n=4$.

(iii) We proceed by recurrence over $n$. In a similar way as in (ii),
we get the result for $n=5$.
Now, we assume that it holds for some $n\geq5$. Again, the same
arguments as in the proof of Corollary \ref{main2} together with
recurrence assumption show that  $\HH^1_{\rm
diff} (\cK(n)^i ; \mathbb{D}^{n+1}_{\lambda,\mu})=0$ if $2(\mu-\lambda)\neq-1,\,0,\,1.$
So, we consider only the cases where $2(\mu-\lambda)=-1,\,0,\,1$, we proceed as in  (i) and we
get the result for $n+1$.
\hfill $\Box$

\medskip

\noindent {\bf Acknowledgements} We are pleased to  thank Dimitry Leites and Valentin Ovsienko for many
stimulating discussions and valuable correspondences.


\end{document}